\documentclass [a4paper,twoside,12pt]{article}

\usepackage[utf8]{inputenc}
\usepackage{amsmath,amsfonts,amssymb}
\usepackage{vmargin,graphicx,theorem}
\usepackage[english]{babel}
\usepackage{enumerate}
\usepackage{color}
\usepackage{pst-fill,pst-grad,pst-plot,pst-eucl,pstricks-add,pst-node}

\setpapersize[portrait]{A4}
\setmarginsrb{1.5cm}{1cm}{1.5cm}{2.5cm}{1.5cm}{0cm}{0.5cm}{2cm}

\newcommand{\disp}{\displaystyle}

\newcommand{\dR}{\ensuremath{\mathbb{R}}}

\newtheorem{ethm}{Theorem}[section]

\newtheorem{ecor}[ethm]{Corollary}

\newtheorem{eprop}[ethm]{Proposition}

\newtheorem{elem}[ethm]{Lemma}

\newtheorem{erem}[ethm]{Remark}

\newcommand{\proofend}{~$\rhd$}
\newcommand{\proofbegin}{~$\lhd$}

\newenvironment{eproof}
               {\noindent {\emph{\textbf{Proof}}}\\\proofbegin~}
               {\proofend\\}

\newcommand{\p}[4]{{#3}\!\left#1{#4}\right#2}

\newcommand{\PAR}[1]{\ensuremath{{\left(#1\right)}}} 

\renewcommand{\phi}{\varphi}

\renewcommand{\geq}{\geqslant}

\newcommand{\entf}[1]{{\rm{Ent}}_{#1}}

\newcommand{\ent}[2]{\p(){\entf{#1}}{#2}}

\def\disp{\displaystyle}

\newcommand{\R}{\dR}

\newcommand{\beq}{\begin{equation}}\newcommand{\eeq}{\end{equation}}
\begin{document}

\title{Dimensional contraction in Wasserstein distance for diffusion semigroups  on a Riemannian manifold}
\author{ Ivan Gentil\thanks{Institut Camille Jordan, Umr Cnrs 5208, Universit\'e Claude Bernard Lyon 1, 43 boulevard du 11 novembre 1918, F-69622 Villeurbanne cedex. gentil@math.univ-lyon1.fr}}

\date{\today}

\maketitle

\abstract{We prove a refined  contraction inequality for diffusion semigroups  with respect to the Wasserstein distance on a compact Riemannian manifold taking account of the dimension.  The result generalizes in a Riemannian context, the dimensional contraction established in~\cite{bgg13} for the Euclidean heat equation. It  is proved by using  a dimensional coercive estimate for the Hodge-de Rham semigroup on 1-forms.   }

\bigskip

\medskip


\section{Introduction}

The von Renesse-Sturm Theorem (c.f. \cite{sturm-vonrenesse})  insures that the contraction of the heat equation on a Riemannian manifold with respect to the Wasserstein distance is equivalent to a uniform lower bound of the Ricci curvature. This result is one of the first equivalence  theorems relating the Wasserstein distance and the Ricci curvature. Actually, there are many  extensions of this result  including the case of a heat equation on a metric measure space and we would like to take account of the dimension into such contraction. 
 
 \medskip

Let us explain the contraction inequality and its extensions with more details. For simplicity, we focus on the heat equation on a Riemannian manifold but all of these results have been proved for a general diffusion semigroup on a Riemannian manifold (or more general spaces). Let  $\Delta$ be the Laplace-Beltrami operator on a smooth Riemannian manifold $(M,g)$ and let $P_t f$ be  the solution of the heat equation $\partial_tu=\Delta u$ with $f$ as  the initial condition. 

  The von Renesse-Sturm Theorem states that:  Let $R\in\R$, the following assertions  are equivalent   (the Wasserstein distance is denoted $W_2$), 
\begin{enumerate}[(i)] 
\item  for any $f,g$ probability densities with respect to the Riemannian measure $dx$ 
$$
 W_2^2(P_tfdx,P_tgdx)\leq e^{-2Rt}W_2^2(fdx,gdx),\,\,\,\forall t\geq0,
 $$
 \item $Ricci_g\geq R$  (uniformly in $M$) where $Ricci_g$ is the Ricci tensor of $(M,g)$. The inequality has to be understood in the sense of inequality between  symmetric tensors. 
 \end{enumerate}

 There are  many proofs and extensions of this result, one can see for instance~\cite{wang-book,otto05,ambrosio-gigli-savare,kuwada10,wang11, savare13,gko13,bglbook,bgl2}.
 
Recently, many extensions have been given taking account of the dimension of the manifold. For instance in~\cite{kuwada,bgl2} the authors prove that if $M$ is a $n$-dimensional Riemannian manifold  with a non-negative Ricci curvature, then for any $f,g$ probability densities with respect to $dx$, and any $s,t\geq0$,  
\begin{equation}
\label{eq-cas-simple}
 W_2^2(P_sfdx,P_tgdx)\leq W_2^2(fdx,gdx)+2n(\sqrt{s}-\sqrt{t})^2,\,\,\, \forall s,t\geq0,
\end{equation}
Non-negative curvature condition has been removed  in~\cite{kuwada,EKS13}. If $M$ is a $n$-dimensional Riemannian manifold, the main extensions are the following:

\begin{itemize}
\item In~\cite{kuwada}, K. Kuwada  proves that the Ricci curvature is bounded from below by $R\in\R$  if and only if for every $s,t\geq0$, 
\begin{equation}
\label{last}
W_2^2(P_t fdx,P_sgdx)\leq A(s,t,R) W_2^2( fdx,gdx)+B(s,t,n,R),\,\,\,\forall s,t\geq0,
\end{equation}
for probability densities $f,g$ with respect to $dx$, for appropriate functions $A,B\geq0$. 

\item In~\cite{EKS13},  M. Erbar, K. Kuwada and K.-T. Sturm prove that the Ricci curvature of the $n$-dimensional manifold $M$ is bounded from below by a constant $R\in\R$ if and only if
\begin{multline}
\label{eq-eks}
s_{\frac Rn}\left(\frac 12 W_2(P_tf dx,P_sg dx)\right)^2
\leq e^{-R(t+s)}\,s_{\frac Rn}\left(\frac 12 W_2(fdx,gdx)\right)^2\\
+\frac nR(1-e^{-R(s+t)})\frac{(\sqrt{t}-\sqrt{s})^2}{2(t+s)},
\end{multline}
for any $s,t\geq0$ and any probability densities $f$ and $g$. Here $s_r(x)=\sin(\sqrt{r}x) / \sqrt{r}$ if $r>0$, $s_r(x)=\sinh(\sqrt{-r}x) / \sqrt{-r}$ if $r<0$ and $s_0(x)=x$, hence recovering~\eqref{eq-cas-simple} when $R=0$.

\item In~\cite{bgg13}  we prove that the classical heat equation in the Euclidean space $\R^n$ satisfies, for any $f$, $g$ probability densities with respect to the Lebesgue measure $\lambda$,
\begin{equation}
\label{eq-chaleur}
W_2^2(P_t f\lambda,P_tg\lambda)\leq W_2^2(f\lambda,g\lambda)-\frac{2}{n} \int_0^t\big(\ent{\lambda}{P_uf}-\ent{\lambda}{P_ug}\big)^2du,\,\,\,\forall t\geq0,
\end{equation}
where $\rm{Ent}$ is the Entropy (it will be defined later). 
\item In the same way, again in~\cite{bgg13}, we prove a more general result for a 
$n$-dimensional Riemannian manifold with a Ricci curvature bounded form below by $R$ and for the  Markov transportation distance $T_2$ (a new distance on measures). We obtain for any $f$ and  $g$, 
$$
\!\!\!\!\!\!\!\!\!\!\!\!\!\!T_2^2(P_t fdx,P_t gdx)\leq e^{-2R t}T_2^2(fdx,gdx)-\frac{2}{n}\int_0^t e^{-2R(t-u)}\PAR{\ent{\mu}{P_u g}-\ent{\mu}{P_u f}}^2du,
$$
for any $t\geq0$.
 \end{itemize}
The goal of this paper is to prove the previous inequality for the Wasserstein distance instead the Markov transportation distance, in other words to extend  inequality~\eqref{eq-chaleur} on a $n$-dimensional Riemannian manifold with a lower bound on the Ricci curvature.

The main result of this paper  can be stated as follows, let $d\mu=e^{-\Psi}dx$ be a probability measure with $\Psi$ a smooth  function on $M$ (for the sequel the manifold will be compact).  We denote  by $(P_t)_{t\geq0}$ the Markov semigroup associated to the generator $L=\Delta-\nabla\Psi\cdot\nabla$.  Then under the curvature-dimension condition 
\begin{equation}
\label{eq-conditioncd}
Ricci_g+\rm{Hess}(\Psi)\geq R+\frac{1}{m-n}\nabla\Psi\otimes\nabla \Psi,
\end{equation}
 for some $R\in\R$ and $m\geq n$ (when $m=n$ then $\Psi=0$),  then  for any probability densities   $f,g$ and any $t\geq0$, 
$$
W_2^2(P_t f\mu,P_t g\mu)\leq e^{-2R t}W_2^2(f\mu,g\mu)-\frac{2}{m}\int_0^t e^{-2R(t-u)}\big[\ent{\mu}{P_u g}-\ent{\mu}{P_u f}\big]^2du. 
$$
The main advantage of such inequality with respect to~\eqref{last} and~\eqref{eq-eks}  is to obtain a contraction inequality with the same time $t$ instead two different times $s$ and $t$. Moreover, the  additional term  is given with a minus sign which shows the improvement  given by the dimension. 

\medskip

The method  to get a dimensional contraction is radically different that one one used in~\cite{EKS13}. Here the strategy to prove such inequality is the Benamou-Brenier dynamical formulation (the Eulerian formulation) of the Wasserstein distance associated to a sharp dimensional estimate on the Hodge-de Rham semigroup. The strategy is closed to the one used in~\cite{otto05}.   In~\cite{EKS13} the authors use in force the definition of the heat equation as a gradient flow of the entropy with respect to the Wasserstein distance. 

\bigskip

The paper is organized as follow. First in Section~\ref{sec-setting},   we recall the Riemannian setting and the Wasserstein distance through  the Benamou-Brenier dynamical formulation.  We need to introduce the Hodge-de Rham operator on forms and its associated semigroup.  In Section~\ref{sec-BLW}, we improve the  Bochner-Lichnerowicz-Weitzenb\"ock identity for 1-forms to get a coercive inequality for the Hodge-de Rham semigroup. Finally in Section~\ref{sec-proof} the main theorem is proved.  

\medskip

For simplicity reasons,  the result will be stated and proved in the context of a compact Riemannian manifold but its generalization  in a metric space,  including the equivalence with respect to the condition~\eqref{eq-conditioncd}, is actually a project with F. Bolley, A. Guillin and K. Kuwada. Moreover, the main theorem  is written for a reversible semigroup but the proof can be adapted to the non-reversible case. In that case the solution of the heat equation is not a gradient flow of the entropy with respect to the Wasserstein distance and the method proposed in~\cite{EKS13} can not be applied. 

\section{Framework and main result}
\label{sec-setting}
\subsection{Geometrics tools}

\noindent{\bf Conventions and notations.} Let $(M,g)$ be a $n$-dimensional, connected, compact and differentiable Riemannian manifold without boundary. We assume for simplicity that  the manifold is $\mathcal C^\infty$. For each $x\in M$ we denote by $T_xM$ the tangent space to $M$ and by $TM$ the whole tangent bundle of $M$.  Moreover  $g_x$ is a symmetric definite positive  quadratic form on $T_xM$, in a local basis of $T_xM$ $(e_i)_{1\leq i\leq n}$, $g_x=(g_{ij})$. (For simplicity, the $x$ dependance of the metric $g$ is omitted.)   For every $x,y\in M$, $d(x,y)$ denotes the usual (Riemannian) distance and $dx$ its measure. In the sequel we will use the Einstein convention of summation over repeated indices: for instance $x_i y^i=\sum_{i=1}^nx_i y^i$, $x^ig_{ij} y^j=\sum_{i,j=1}^nx^ig_{ij} y^j$.

The Riemannian scalar product between two tensors $X$ and $Y$ is denoted $X\cdot Y$, its  associated norm is noted $|X|$ (depending on $g$). For instance, locally in a basis $(e_i)_{1\leq i\leq n}$ for two smooth functions $f,h:M\mapsto \R$, $\nabla f\cdot\nabla h=\partial_if g^{ij}\partial _jh$ where $(g^{ij})=(g_{ij})^{-1}$.    As usual the covariant derivative in the direction $e_i$ of a tensor $X$ (vector field or form) is noted $\nabla_iX$. The geometric musicology will be used in force, for instance $\nabla^i f=g^{ij}\nabla_j f=g^{ij}\partial_j f$ or $\nabla^iX=g^{ij}\nabla_j X$. If $\omega$ is a 1-form then $\omega^*$ is its dual representation as a vector field with components $\omega^i=g^{ij}\omega_j$. 
The Ricci tensor of  $(M,g)$ is noted $Ricci_g$. 

The Laplace-Beltrami operator $\Delta$ is acting on smooth functions $f$,  
$$
\Delta f=\nabla\cdot \nabla f,
$$
where $\nabla\cdot$ is the divergence operator acting on vector fields, $\nabla\cdot X=\nabla_i X^i$ for every vector field $X$. (We use the analyst's convention with respect to the sign.) A smooth function (or form) is a $\mathcal C^\infty$-function (or form).  The divergence operator $\nabla\cdot$ satisfies for any smooth vector fields  $X$ and any smooth functions $f$, 
$$
\int_M  \nabla \cdot X fdx=-\int_M  \nabla f\cdot X  dx.
$$ 
The Laplace-Betrami operator can also be written as $\Delta f=\delta df$ 
where $\delta$ is the divergence operator on 1-form
$$
\delta \omega =\nabla_i \omega^i=g^{ij}\nabla_i \omega_j=g^{ij}(\partial_i\omega_j-\Gamma_{ij}^p \omega_p), 
$$
where $\Gamma_{ij}^p$ are Christoffel symbols. Operators $\delta$ and $\nabla\cdot$ are related by the formula  $\delta \omega=\nabla\cdot \omega^*$, for any 1-forms $\omega$. 

\medskip
\noindent{\bf The Markov (or heat) semigroup.} Let $\Psi:M\mapsto \R$ be a  fixed $\mathcal C^\infty$-function. Let $\mu(dx)=e^{-\Psi}dx$ and since $M$ is compact  we can  assume that $\mu$ is a probability measure.  Let $Lf=\Delta f-\nabla\Psi\cdot\nabla f$, for any smooth functions $f$. Since the manifold is compact, the operator $L$ defines a unique Markov semigroup $(P_t)_{t\geq0}$ on $L^2(\mu)$ and it is called a Markov generator.    This semigroup is symmetric in $L^2(\mu)$ and  for any  $f$, $P_tf$ is a solution of the equation $\partial_t u=L u$ with $f$ as the initial condition. If $f$ is a probability density  (with respect to  $\mu$) then for all $t\geq0$, $P_t f \mu$ remains a probability measure.  Finally, for any smooth  functions $f$ and $g$ on $M$, 
$$
\int_M   fL gd\mu=\int_M  g L f\mu=-\int_M  \nabla f\cdot\nabla g d\mu=-\int_M \Gamma(f,g)d\mu,
$$
where $\Gamma$ is  Carr\'e du champ operator defined on functions, $\Gamma(f,g)=\nabla f\cdot\nabla g$.  If $X$ is a vector field, we define $\nabla_\Psi\cdot X=\nabla\cdot X-\nabla \Psi\cdot X$, and the generator $L$ takes then  the form 
$$
L=\nabla_\Psi \cdot\nabla. 
$$
As before we note by $\delta_\Psi$ the divergence operator acting on forms: 
\begin{equation}
\label{eq-def-divpsi}
\delta_\Psi \omega=\nabla_\Psi\cdot \omega^*.
\end{equation}
It satisfies the integration by parts formula, for any smooth 1-forms $\omega$ and functions $f$,
$$
\int_M\delta_\Psi \omega f\,d\mu=-\int_M \omega\cdot d fd\mu, 
$$
where $\omega\cdot d f=\omega_idf^i$ (the inner product between the two 1-forms). 
 
The triple $(M,\mu,\Gamma)$ is a compact Markov triple as defined in~\cite[Chap. 3]{bglbook}.

\medskip
\noindent{\bf The Hodge-de Rham semigroup.} Connected to the Markov semigroup $(P_t)$  (associated to the generator $L$)  one can define the Hodge-de Rham semigroup. As explained in this context  in~\cite{bakry87} the  (modified)  Hodge-de Rham operator is acting on smooth 1-forms, 
\begin{equation}
\label{eq-defhodge}
\overrightarrow{L}=(d^{(0)}\delta_\Psi^{(1)} +\delta_\Psi^{(2)}d^{(1)}), 
\end{equation}
where $d^{(i)}$ is the differential operator acting on $i$-forms and $\delta_\Psi^{(i)}$ is its adjoint operator in $L^2(\mu)$ with respect to the usual inner product on $i$-form : $\int d^{(i)}\omega\cdot\eta \,d\mu=\int\omega\cdot\delta_\Psi^{(i)}\eta\, d\mu$ for any $i$-forms $\omega$ and $(i+1)$-forms $\eta$.  In the sequel,  we omit the exponent $(i)$. For computations,  we  use the Hodge-de Rham operator toward the Weitzenb\"ock formula, that is for any $1\leq i\leq n$, 
\begin{equation}
\label{eq-weitzenbock}
(\overrightarrow{L}\omega)_i=\nabla^k\nabla_k\omega_i-\big(\nabla_{\nabla\Psi}\omega\big)_i- Ricci(L)(\omega^*,e_i)=\nabla^k\nabla_k\omega_i-\nabla^k\Psi\nabla_{k}\omega_i- Ricci(L)(\omega^*,e_i),
\end{equation}
where $Ricci(L)=Ricci_g+\rm{Hess}(\Psi)$ is the so-called Bakry-\'Emery tensor (see for instance~\cite[Prop. 1.5]{bakry87}). Again, we use the analyst's convention with respect to the sign. If $\Psi=0$, then  $L=\Delta$ and $\overrightarrow{L}$ is the usual Hodge-de Rham operator noted $\overrightarrow{\Delta}$.

Since $M$ is compact,  the operator $\overrightarrow{L}$ induces a  semigroup $(R_t)_{t\geq0}$ on 1-forms. It is also symmetric in ${L}^2(d\mu)$,  for any smooth 1-forms $\omega$ and $\eta$,
$$
\int_M  \omega\cdot\overrightarrow{L}\eta\, d\mu=\int_M  \eta\cdot\overrightarrow{L}\omega\, d\mu.
$$
Then for any smooth 1-forms $\omega$, $R_t\omega$ is the solution of  $\partial_t {u}=\overrightarrow{L}{u}$ where ${u}:[0,\infty)\times M\mapsto TM^*$ ($TM^*$ is the cotangent bundle) with  $\omega$ as the initial condition.  The details of the construction of the Hodge-de Rham semigroup can be found in~\cite{bakry87} (see also the references therein).

The Hodge-de Rham semigroup is related to the Markov generator  by the following commutation property: for any smooth 1-forms $\omega$ and $t\geq0$, 
\begin{equation}
\label{eq-commutation}
P_t \delta_\Psi\omega=\delta_\Psi R_t \omega. 
\end{equation}
The easiest way to prove this fundamental identity is to use the definition~\eqref{eq-defhodge} and the identity  $\delta_\Psi^{(1)}\delta_\Psi^{(2)}=0$.

\subsection{The Wasserstein distance}
\label{sec-wasser}
Let $\mathcal P (M)$ be the set of probability measures in $M$. The Wasserstein distance between two  probability measures $\nu_1,\nu_2\in\mathcal P(M)$,   is defined by 
$$
W_2(\nu_1,\nu_2)=\inf\Big(\int_{M\times M} d^2(x,y)d\pi(x,y)\Big)^{1/2},
$$
where the infimum runs over all probability measures $\pi$ in $M\times M$ with marginals $\nu_1$ and $\nu_2$. We refer to the monumental work~\cite{villani-book1} for a reference presentation of this distance, its interplay with the optimal transportation problem and many other issues. 

The Wasserstein distance has a dynamical formulation:  for any  probabilities measure $\nu_1,\nu_2\in\mathcal P(M)$, 
$$
W_2^2(\nu_1,\nu_2)=\inf \int_0^1\int_M |\eta_s|^2d\mu_s ds, 
$$
where the infimum is running over all paths of probabilities  $(\mu_s)_{s\in[0,1]}$ and $\eta_s\in TM^*$ satisfying in the distributional sense 
$$
\left\{
\begin{array}{l}
\partial_s\mu_s+\delta (\mu_s\eta_s)=0\\
\mu_0=\nu_,\,\,\,\mu_1=\nu_2.
\end{array}
\right.
$$
The Euclidean case has been proved by Benamou-Brenier in~\cite{benamoubrenier} and the Riemannian case by F. Otto and M. Westdickenberg an the other hand by L. De Pascale, M. S. Gelli and L. Granieri in~\cite{otto05,pascalegelli}. 

Let $\omega_s=\eta_s \rho_s$ where $\frac{d\mu_s}{d\mu}=\rho_s$,  in~\cite{otto05} the authors state that for any $f$ and $g$, smooth  probability densities (with respect to  $\mu$), 
\begin{equation}
\label{eq-bbotto}
W_2^2(f\mu,g\mu)=\inf \int_0^1\int_M\frac{|\omega_s|^2}{\rho_s}d\mu ds, 
\end{equation}
where the infimum is running over smooth couples  $(\rho_s,\omega_s)_{s\in[0,1]}$ where for any $s\in[0,1]$, $\rho_s$ is a positive probability density (with respect to $\mu$) and $\omega_s$ is a 1-form,  satisfying  
$$
\left\{
\begin{array}{l}
\partial_s\rho_s+\delta_\Psi \omega_s=0\\
\rho_0=f,\,\,\,
\rho_1=g.
\end{array}
\right.
$$
This dynamical formulation of optimal transportation has been manly used to get contraction result or Evolutional Variational Inequality (see e.g. \cite{dns1,dns2,bgg13}). 

\subsection{Main result}

We can now state the main result of this paper. 
\begin{ethm}[Dimensional contraction in Wasserstein distance]
\label{theo-main}
Let $(M,g)$ be a   $\mathcal  C^\infty$, $n$-dimensional, connected and compact  Riemannian manifold and let $\Psi:M\mapsto \R$ be a $\mathcal C^\infty$-function.  We assume that there exits $R\in\R$ and $m\geq n$  such that uniformly in $M$ (if $m=n$ then we impose that $\Psi=0$), 
\begin{equation}
\label{eq-prems}
Ricci(L)\geq R+\frac{1}{m-n}\nabla\Psi\otimes\nabla \Psi.
\end{equation}
Then for any  smooth probability densities  $f,g$ with respect to $\mu$ and any $t\geq0$, 
\begin{equation}
\label{eq-maintheo}
W_2^2(P_t f\mu,P_t g\mu)\leq e^{-2R t}W_2^2(f\mu,g\mu)
-\frac{2}{m}\int_0^t e^{-2R(t-u)}\big[\ent{\mu}{P_u g}-\ent{\mu}{P_u f}\big]^2du,
\end{equation}
where  $\ent{\mu}{h}=\int_M  h\log h d\mu$ for every probability density $h$ (with respect to $\mu$). 

\end{ethm}

\begin{erem}
When $m\rightarrow\infty$ we recover the von Renesse-Sturm result, the exponential contraction of the heat equation with respect to the Wasserstein distance. When $\Psi=0$, one can choose  $m=n$, and then the Laplace-Beltrami operator satisfies the condition~\eqref{eq-prems} under a lower bound on the Ricci curvature.

As it will be explained in Remark~\ref{rem-be}, condition~\eqref{eq-prems} is equivalent to the so-called Bakry-\'Emery curvature-dimension  condition $CD(R,m)$. 
\end{erem}

\section{Coercive inequality for the Hodge-de Rham semigroup}
\label{sec-BLW}

The next proposition state a refined Bochner-Lichnerowicz-Weitzenb\"ock formula for 1-forms. 

\begin{eprop}[Refined Bochner-Lichnerowicz-Weitzenb\"ock formula]
\label{prop-blwgene}
For any smooth 1-forms $\alpha,\eta$ and any $b\in\R$, 
\begin{equation}
\label{eq-supergene}
L\frac{|\eta|^2}{2}-\eta\cdot\overrightarrow{L}{\eta}+2b\,\alpha\cdot d|\eta|^2+4b^2\,|\alpha|^2|\eta|^2=|\nabla\eta+2b\,\alpha\otimes\eta|^2+Ricci(L)(\eta^*,\eta^*),
\end{equation}
where $|\nabla\eta+2b\,\alpha\otimes\eta|$ has to be understood as the norm of the  2-tensor $\nabla_i\eta_j+2b\,\alpha_i\eta_j$.
\end{eprop}

When $\Psi=0$, $\overrightarrow{\Delta}$ be the usual Hodge-de Rham operator. The Weitzenb\"ock~\eqref{eq-weitzenbock} implies
\begin{equation}
\label{eq-aidew}
\overrightarrow{L}\omega=\overrightarrow{\Delta}\omega-\nabla_{\nabla\Psi}\omega-{\rm Hess}(\Psi)(\omega^*,\cdot).
\end{equation}
According to Lemma~\ref{lemm0}, equation~\eqref{eq-supergene} is equivalent to the  identity  
$$
\Delta\frac{|\eta|^2}{2}-\eta\cdot\overrightarrow{\Delta}{\eta}+2b\,\alpha\cdot d|\eta|^2+4b^2\,|\alpha|^2|\eta|^2=|\nabla\eta+2b\,\alpha\otimes\eta|^2+Ricci_g(\eta^*,\eta^*),
$$
which will be proved in  Proposition~\ref{prop-blw}. The main difficulty is to prove it for all $b\in\R$ since when $b=0$, equation~\eqref{eq-supergene} is classic.

Next Lemma can be found for instance in~\cite[p. 3]{lichnerowicz58}. 
\begin{elem}
\label{lemm0}
For any 1-forms $\eta$ and any $1\leq k\leq n$, $\Big(d\frac{|\eta|^2}{2}\Big)_k=\eta^i\nabla_k \eta_i$.
\end{elem}

\begin{eprop}
\label{prop-blw}
For any smooth 1-forms $\alpha$, $\eta$ and any $b\in\R$, 
\begin{equation}
\label{eq-super}
\Delta\frac{|\eta|^2}{2}-\eta\cdot\overrightarrow{\Delta}{\eta}+2b\,\alpha\cdot d|\eta|^2+4b^2\,|\alpha|^2|\eta|^2=|\nabla\eta+2b\,\alpha\otimes\eta|^2+Ricci_g(\eta^*,\eta^*).
\end{equation}
\end{eprop}
\begin{eproof}
The Bochner-Lichnerowicz-Weitzenb\"ock formula (c.f.~\cite[P. 3]{lichnerowicz58}) insures that for any smooth  1-forms $\omega$, 
\begin{equation}
\label{eq-blw}
\Delta \frac{|\omega |^2}{2}-\omega\cdot\overrightarrow{\Delta}\omega=|\nabla\omega |^2+Ricci_g(\omega^*,\omega^*),
\end{equation}
recall that  $|\nabla\omega |^2={\nabla_i\omega_j g^{il}g^{jk}\nabla_l\omega_k}$ is the square of the norm of the 2-tensor $\nabla_i\omega_j$.

The idea is to change variables into this formula. We would like to prove this result at some $x_0\in M$ which is supposed to be fixed. Let
$$
\omega=(bg+1)\eta+b\,f\alpha,
$$
where   $\eta$ and $\alpha$ are 1-forms and $f,g$ are actually smooth functions satisfying $f(x_0)=g(x_0)=0$. 

First we have 
$$
|\omega|^2=(bg+1)^2|\eta|^2+b^2\,f^2|\alpha|^2+2b(bg+1)f\alpha\cdot \eta.
$$
Since, for any smooth functions $F$ and $G$ we have $\Delta(FG)=2\,\Gamma(F,G)+F\Delta G+G\Delta F$,
where  $\Gamma$ is the carr\'e du champ operator,  a straight forward computation gives at the point $x_0$, 
\begin{multline*}
\Delta |\omega |^2=4b\Gamma( g,|\eta|^2)+2b|\eta|^2\Delta g+2b^2|\eta|^2|\nabla g|^2+\Delta |\eta|^2+2b^2|\alpha|^2|\nabla f|^2\\
+4b\Gamma(f, \eta\cdot\alpha)+4b^2(\eta\cdot\alpha) \Gamma(f,g)+2b(\eta\cdot\alpha)\, \Delta f.
\end{multline*}
Since $\omega(x_0)=\eta(x_0)$, we have in $x_0$, 
$$
2\omega\cdot\overrightarrow{\Delta}\omega=2b\eta\overrightarrow{\Delta}(g\eta)+2\eta\overrightarrow{\Delta}{\eta}+2b\eta\overrightarrow{\Delta}(f\alpha).
$$
Thanks to Lemma~\ref{lem-diffusion}, at the point $x_0$ (recall that $f(x_0)=g(x_0)=0$),
$$
2\omega\cdot\overrightarrow{\Delta}\omega=2b|\eta|^2{\Delta}g+4b\eta\cdot\nabla_{\nabla g}\eta+2\eta\overrightarrow{\Delta}{\eta}+2b\eta\cdot\alpha\Delta f+4b\eta\cdot\nabla_{\nabla f}\alpha.
$$
Moreover,  at the point $x_0$,
$$
|\nabla\omega|^2=b^2|\eta|^2\Gamma(g)+|\nabla\eta|^2+b^2\Gamma(f)|\alpha|^2+2b\eta\cdot \nabla_{\nabla g} \eta
+2b\alpha\cdot \nabla_{\nabla f}\eta+2b^2\eta\cdot \alpha\,\Gamma(f,g),
$$
and  $Ricci_g(\omega^*,\omega^*)=Ricci_g(\eta^*,\eta^*)$. At the point $x_0$, the identity~\eqref{eq-blw} applied to $\omega$, becomes 
\begin{multline*}
\Delta \frac{|\eta|^2}{2}-\eta\overrightarrow{\Delta}{\eta}+2b\Gamma(g,|\eta|^2)+b^2|\eta|^2\Gamma(g)+b^2|\alpha|^2\Gamma(f)
+2b\Gamma(f,\eta\cdot\alpha)\\
+2b^2(\eta\cdot\alpha)\Gamma ( f, g)
-2b\eta\cdot\nabla_{\nabla g}\eta-2b\,\eta\cdot\nabla_{\nabla f}\alpha
=b^2|\eta|^2\Gamma(g)+|\nabla \eta|^2\\
+b^2\Gamma(f)|\alpha|^2
+2b\eta\cdot\nabla_{\nabla g}\eta
+2b\alpha\cdot\nabla_{\nabla f}\eta+2b^2\eta\cdot\alpha \Gamma(g,f)+Ricci_g(\eta^*,\eta^*)
\end{multline*}
Let now assume that  $f$ and $g$ satisfied moreover $df(x_0)=\eta(x_0)$ and $dg(x_0)=\alpha(x_0)$ (which is always possible),  we obtain
\begin{multline*}
\Delta \frac{|\eta|^2}{2}-\eta\overrightarrow{\Delta}{\eta}+2b\,\alpha\cdot\nabla|\eta|^2+4b^2|\eta|^2|\alpha|^2
=4b^2\,|\eta|^2|\alpha|^2-2b\,\eta\cdot \nabla (\eta\cdot\alpha)
+2b\,\eta\cdot\nabla_{\alpha}\eta\\+2b\,\eta\cdot\nabla_{\eta}\alpha+|\nabla \eta|^2
+2b\,\eta\cdot\nabla_\alpha\eta
+2b\,\alpha\cdot\nabla_\eta\eta+Ricci_g(\eta^*,\eta^*).
\end{multline*}
Then Lemma~\ref{lem-classique} implies that 
\begin{multline*}
4b^2\,|\eta|^2|\alpha|^2-2b\,\eta\cdot d (\eta\cdot\alpha)
+2b\,\eta\cdot\nabla_{\alpha}\eta+2b\,\eta\cdot\nabla_{\eta}\alpha+
|\nabla\eta|^2
+2b\,\eta\cdot\nabla_\alpha\eta
+2b\,\alpha\cdot \nabla_\eta\eta\\
=4b^2\,|\eta|^2|\alpha|^2+4b\,\eta\cdot\nabla_\alpha\eta+|\nabla \eta|^2=|\nabla\eta+2b\,\alpha\otimes\eta|^2,
\end{multline*}
which implies the identity~\eqref{eq-super}.
\end{eproof}
The next lemma is classic  in the Riemannian context and we skip the proof.
\begin{elem}
\label{lem-classique}
For any 1-forms $\eta$ and $\alpha$, $d(\eta\cdot\alpha)=\nabla\eta\cdot\alpha+\nabla\alpha\cdot\eta$, e.g. in coordinates 
$$
(d(\eta\cdot\alpha))_i=\alpha^j\nabla_i\eta_j+\eta^j\nabla_i\alpha_j.
$$
\end{elem}

The next result can be seen as a kind of diffusion property as defined in~\cite{bglbook}. The result is probably classic but I didn't find it in the literature. 
\begin{elem}
\label{lem-diffusion}
For any smooth functions $f$ and 1-form $\omega$, 
\begin{equation}
\label{eq-diffusion}
\overrightarrow{\Delta}(f\omega)=f\overrightarrow{\Delta}(\omega)+\omega\Delta f+2\nabla_{\nabla f}\eta.
\end{equation}
In other words, for any $i$, $
\big(\overrightarrow{\Delta}(f\omega)\big)_i=f\big(\overrightarrow{\Delta}(\omega)\big)_i+\omega_i\Delta f+2\nabla ^jf\nabla_{j}\eta_i$.
\end{elem}

\begin{eproof}
The Weitzenb\"ock formula~\eqref{eq-weitzenbock} insures that  
$$
\big(\overrightarrow{\Delta}(f\omega)\big)_i=\nabla^k\nabla_k(f\omega)_i- Ricci_g(f\omega^*,e_i)=\nabla^k\nabla_k(f\omega)_i- fRicci_g(\omega^*,e_i).
$$
Now, using in force the formula $\nabla_k(f\omega)_i=\omega_i\partial_kf+f\nabla_k\omega_i$, we get
$$
\begin{array}{rl}
\nabla^k\nabla_k(f\omega)_i=&g^{kl}\nabla_l\nabla_k(f\omega)_i\\
=&\disp g^{kl}(\nabla_k\nabla_l f)\omega_i+fg^{kl}\nabla_l\nabla_k\omega_i+2g^{kl} \partial_kf\nabla_l\omega_i\\
=&\disp\omega_i\Delta f+f\nabla^j\nabla_j\omega_i+2\nabla^j f\nabla_j\omega_i=\disp\omega_i\Delta f+f\nabla^j\nabla_j\omega_i+2\nabla_{\nabla f}\omega_i,
\end{array}
$$
which implies~\eqref{eq-diffusion}.
\end{eproof}
Finally we can state our main estimate.
\begin{ecor}
\label{cor-blw}
Assume that there exits $R\in\R$ and $m\geq n$  such that uniformly in  $M$  (if $m=n$ then we impose that $\Psi=0$),
\begin{equation}
\label{eq-dis}
Ricci(L)\geq R+\frac{1}{m-n}\nabla\Psi\otimes\nabla \Psi.
\end{equation}
Then for any smooth 1-forms $\eta$, $\alpha$ and $b\in\R$, (with $\delta_\Psi$ defined in~\eqref{eq-def-divpsi}),

\begin{equation}
\label{eq-corcd}
L\frac{|\eta|^2}{2}-\eta\cdot\overrightarrow{L}{\eta}+2b\alpha\cdot d|\eta|^2+4b^2|\alpha|^2|\eta|^2\geq\frac{1}{m}(\delta_\Psi \eta+2b\alpha\cdot\eta)^2+R|\eta|^2.
\end{equation}
\end{ecor}
\begin{eproof}
From~\eqref{eq-supergene}  we only have to verify 
$$
|\nabla\eta+2b\alpha\otimes\eta|^2+Ricci(L)(\eta^*,\eta^*)\geq\frac{1}{m}(\delta_\Psi \eta+2b\alpha\cdot\eta)^2+R|\eta|^2. 
$$
Let $x\in M$ and let assume that $(e_i)_{1\leq i\leq n}$ is an orthonormal basis of $T_xM$. Then the left hand side can be written 
$$
|\nabla\eta+2b\alpha\otimes\eta|^2=\sum_{i,j=1}^n(\nabla_i\eta_j+2b\alpha_i\eta_j)^2
\geq \frac1n\big(\sum_{i=1}^n\nabla_i\eta_i+2b\alpha\cdot\eta\big)^2,
$$
from the Cauchy-Schwartz inequality. So it remains to prove that for any $t\in\R$
$$
\frac1n(t+2b\alpha\cdot\eta)^2+Ricci(L)(\eta^*,\eta^*)\geq\frac{1}{m}(t-d\Psi\cdot \eta+2b\,\alpha\cdot\eta)^2+R|\eta|^2,
$$
where $t=\delta\eta=\sum_{i=1}^n\nabla_i\eta_i$ (recall that $(e_i)$ is an orthonormal basis). The formula is valid since $m\geq n$ and from~\eqref{eq-dis} the discriminant of this two order polynomial function in variable $t$ is non-positive (it doesn't depend on the parameter $b$). 

As proposed by the referee the last inequality can be stated directly by using the quadratic inequality 
$$
\frac{1}{m-n}x^2+\frac{1}{n}{y^2}\geq \frac{1}{m}(x+y)^2, 
$$
with some reals $x$ and $y$. 
\end{eproof}

\begin{erem}[Link with the curvature-dimension condition]
\label{rem-be}
 The so-called Bakry-\'Emery curvature-dimension condition $CD(R,m)$ for an operator $L$ is satisfied when for every smooth function $f$,
$$
\Gamma_2(f)\geq R\,\Gamma(f)+\frac1m(Lf)^2,
$$ 
where 
$$
\Gamma_2(f)=\Gamma_{2}(f,f)=\frac12\big(L\Gamma(f)-2\Gamma(f,Lf)\big).
$$
The same procedure in the case of closed 1-forms has been stated in~\cite{bakrysaintflour} (see also~\cite[Chap. 5]{logsob}) in the context of $\Gamma_2$-calculus. It is proved that if $\eta=df$ and $\alpha=dg$,  then under the curvature-dimension inequality $CD(R,m)$, and any $b\in\R$, 
\begin{equation}
\label{eq-gamma2}
\Gamma_2(f)+2b\Gamma(f,\Gamma(g))+4b^2\Gamma(f)\Gamma(g)\geq\frac{1}{n}(Lf+2b\Gamma(f,g))^2+\rho\Gamma(f),
\end{equation}
which is~\eqref{eq-corcd} for closed 1-form. 
Moreover, inequality~\eqref{eq-gamma2} for every function $f$ (with $b=0$) is equivalent to $CD(R,m)$. Since inequality~\eqref{eq-corcd} is a generalization of~\eqref{eq-gamma2}, it is also equivalent to~\eqref{eq-dis} and then to $CD(R,m)$. 
\end{erem}
We can now give the main estimation of our semigroups: 
\begin{ethm}[Coercive estimation]
\label{thm-coercive}
Assume that there exits $R\in\R$, $m\geq n$  such that (if $m=n$ then we impose that $\Psi=0$),
\begin{equation}
\label{eq-disbis}
Ricci(L)\geq R+\frac{1}{m-n}\nabla\Psi\otimes\nabla \Psi.
\end{equation}
Then for any smooth 1-forms $\omega$ and smooth  functions $g>0$, for any $t\geq0$, 
\begin{equation}
\label{eq-blwint}
\frac{|R_t\omega|^2}{P_t g}\leq e^{-2R t}P_t\Big(\frac{|\omega|^2}{g}\Big)-\frac{2}{m}\int_0^t \frac{e^{-2Ru}}{P_tg}\big[ P_t\delta_\Psi \omega -P_u(d (\log P_{t-u} g)\cdot R_{t-u}\omega)\big]^2 du.
\end{equation}
\end{ethm}

\begin{eproof}
Since $M$ is compact,  one can assume that there exists $\varepsilon>0$ such that $g=f+\epsilon$ with $f>0$ and then $\varepsilon\rightarrow0$ in~\eqref{eq-blwint}. Let $\omega$ be a smooth 1-form and  $t\geq0$. For any $s\in[0,t]$, we define 
$$
\Lambda(s)=P_s\Big(\frac{|R_{t-s}\omega|^2}{P_{t-s} g}\Big).
$$
For any $s\in[0,1]$, 
$$
\Lambda'(s)=P_s\Big(L\Big(\frac{|R_{t-s}\omega|^2}{P_{t-s} g}\Big)-2\frac{R_{t-s}\omega\cdot\overrightarrow{L}R_{t-s}\omega }{P_{t-s} g}+L P_{t-s} g\frac{|R_{t-s}\omega|^2}{(P_{t-s} g)^2}\Big).
$$
Since for any smooth functions $F,G$, $L(FG)=2\Gamma( F, G)+FL G+GL F$,
the identity becomes, with $\eta=R_{t-s}\omega$ and $G=P_{t-s} g$, 
\begin{multline*}
\Lambda'(s)=P_s\Big(-\frac{2}{G^2}\Gamma( |\eta|^2,G)-\frac{2}{G}\eta\cdot\overrightarrow{L}\eta +\frac{2}{G^3}|\eta|^2\Gamma( G)+\frac1GL|\eta|^2\Big)\\
=P_s\Big[\frac2G\Big(L\frac{|\eta|^2}{2}-\eta\cdot\overrightarrow{L}\eta-\Gamma(|\eta|^2,\log G)+|\eta|^2|\Gamma( \log G)\Big)\Big].
\end{multline*}
Corollary~\ref{cor-blw} appied to $\alpha=d \log G$ and $b=-1/2$, implies 
$$
\Lambda'(s)\geq \frac{2}{m}P_s\Big(\frac{1}{G}\big(\delta_\Psi \eta-d (\log G)\cdot\eta\big)^2\Big)+2R P_s\Big(\frac{|\eta|^2}{G}\Big).
$$
Thanks to the  Cauchy-Schwarz inequality,
 $$
 P_s\Big(\frac{1}{G}\big(\delta_\Psi \eta-d (\log G)\cdot\eta\big)^2\Big)\geq \frac{1}{P_s G}{\Big[P_s\big(\delta_\Psi \eta-d (\log G)\cdot\eta\big)\Big]^2}. 
 $$
 Since $P_s G=P_{s}P_{t-s}g=P_t g$ and from~\eqref{eq-commutation}, $P_s\delta_\Psi R_{t-s}\omega=P_t\delta_\Psi\omega$, the inequality becomes
$$
(\Lambda(s)e^{-2R s})'\geq \frac{2}{m}\frac{e^{-2Rs}}{P_tg}\Big[P_t\delta_\Psi\omega-P_s\big(d (\log P_{t-s}g)\cdot R_{t-s}\omega\big)\Big]^2.
$$
The integration  over $s\in[0,t]$ of the previous inequality  implies~\eqref{eq-blwint}.
\end{eproof}

\section{Proof of Theorem~\ref{theo-main}}
\label{sec-proof}

\noindent
\!\!\proofbegin~{\bf Regularity assumption.} Let  $f_\varepsilon=(P_\varepsilon f+\varepsilon)/(1+\varepsilon)$ and $g_\varepsilon=(P_\varepsilon g+\varepsilon)/(1+\varepsilon)$, for $\varepsilon>0$. The probability measure  $f_\varepsilon\mu$ (resp. $g_\varepsilon\mu$) converges weakly to $f\mu$ (resp. $g\mu$), when $\varepsilon\rightarrow0$  and since $M$ is compact $W_2^2(f_\varepsilon\mu,g_\varepsilon\mu)$ converges to $W_2(f\mu,g\mu)$.  The same is also true for  $W_2^2(P_tf_\varepsilon\mu,P_tg_\varepsilon\mu)$.   Thus,  on can assume that   $f$ and $g$ are two smooth functions satisfying $f,g\geq\varepsilon$ for some $\varepsilon>0$. 

\medskip

\noindent{\bf Contraction inequality.} Let $f,g\geq\varepsilon$, be two smooth functions and let $(\rho_s,\omega_s)_{s\in[0,1]}$ be a smooth couple satisfying 
\begin{equation}
\label{eq-condic2}
\left\{
\begin{array}{l}
\disp\partial_s\rho_s+\delta_\Psi \omega_s=0\\
\rho_0=f,\,\,\,\,\,
\rho_1=g.
\end{array}
\right.
\end{equation}
For any $t\geq0$, from the commutation property~\eqref{eq-commutation},  the couple $(P_t\rho_s,R_t\omega_s)_{s\in[0,1]}$ satisfies 
\begin{equation}
\label{eq-condic}
\left\{
\begin{array}{l}
\partial_sP_t\rho_s+\delta_\Psi R_t\omega_s=0\\
P_t\rho_0=P_tf,\,\,\,\,\,P_t\rho_1=P_tg.
\end{array}
\right.
\end{equation}
On can apply  Theorem~\ref{thm-coercive} to get 
\begin{multline}
\label{eq-avance}
\int_0^1\int_M\frac{|R_t\omega_s|^2}{P_t\rho_s}d\mu ds\leq e^{-2R t}\int_0^1\int_M \frac{|\omega_s|^2}{\rho_s}d\mu ds\\
-\frac{2}{m}\int_0^1\int_0^te^{-2Ru} \int_M \frac{\big[ P_t\delta_\Psi \omega_s -P_u(d (\log P_{t-u} \rho_s)\cdot R_{t-u}\omega_s)\big]^2}{P_t\rho_s} d\mu  du ds.
\end{multline}
Cauchy-Schwarz inequality implies
\begin{multline*}
\int_M \frac{\big[ P_t\delta_\Psi \omega_s -P_u(d (\log P_{t-u} \rho_s)\cdot R_{t-u}\omega_s)\big]^2}{P_t\rho_s} d\mu\\ \geq
\frac{1}{\int P_t\rho_sd\mu}\Big[\int_M\!\!  \Big(P_t\delta_\Psi \omega_s -P_u(d (\log P_{t-u} \rho_s)\cdot R_{t-u}\omega_s)\Big)d\mu\Big]^2
\!\!=\!\!\Big[\int_M\!\!  d (\log P_{t-u} \rho_s)\cdot R_{t-u}\omega_s d\mu \Big]^2,
\end{multline*}
since $\int_M  P_t\delta_\Psi \omega_sd\mu=\int_M  \delta_\Psi \omega_sd\mu=0$. Integrating over $s\in[0,1]$, thanks again to the Cauchy-Schwartz inequality, 
\begin{multline*}
\int_0^1\Big[\int_M  \nabla (\log P_{t-u} \rho_s)\cdot R_{t-u}\omega_s d\mu \Big]^2ds\geq \Big[\int_0^1\int_M  d (\log P_{t-u} \rho_s)\cdot R_{t-u}\omega_s d\mu ds\Big]^2\\
=\big[\ent{\mu}{P_{t-u}f}-\ent{\mu}{P_{t-u}g}\big]^2,
\end{multline*}
since coming from~\eqref{eq-condic} we have (since $M$ is compact, the next integration by parts is valid)
\begin{multline*}
\frac{d}{ds}\ent{\mu}{P_{t-u}\rho_s}=\frac{d}{ds}\int_M  P_{t-u}\rho_s\log P_{t-u}\rho_s d\mu=\int_M  \partial_s P_{t-u}\rho_s\log P_{t-u}\rho_s d\mu\\
=-\int_M  \delta_\Psi R_{t-u}\omega_s\log P_{t-u}\rho_s d\mu=\int_M  R_{t-u}\omega_s\cdot d (\log P_{t-u} \rho_s) d\mu.
\end{multline*}
The inequality~\eqref{eq-avance} becomes 
$$
\int_0^1\int_M\frac{|R_t\omega_s|^2}{P_t\rho_s}d\mu ds\leq e^{-2R t}\int_0^1\int_M \frac{|\omega_s|^2}{\rho_s}d\mu ds-\frac2m\int_0^te^{-2Ru}\big[\ent{\mu}{P_{t-u}f}-\ent{\mu}{P_{t-u}g}\big]^2du
$$
Thanks to the Brenier-Benamou formulation~\eqref{eq-bbotto} with respect  $P_tf\mu$ and $P_tg\mu$ and formula~\eqref{eq-condic}, 
$$
W_2^2(P_tf\mu,P_tg\mu)\leq e^{-2R t}\int_0^1\int_M \frac{|\omega_s|^2}{\rho_s}d\mu ds-\frac2m\int_0^te^{-2Ru}\big[\ent{\mu}{P_{t-u}f}-\ent{\mu}{P_{t-u}g}\big]^2du.
$$
Taking now the infimum over all couples  $(\rho_s,\omega_s)_{s\in[0,1]}$ satisfying~\eqref{eq-condic2},  
$$
W_2^2(P_tf\mu,P_tg\mu)\leq e^{-2R t}W_2^2(f\mu,g\mu)-\frac2m\int_0^te^{-2Ru}\big[\ent{\mu}{P_{t-u}f}-\ent{\mu}{P_{t-u}g}\big]^2du, 
$$
 which is the inequality desired changing $t-u$ by $u$.
\proofend
\bigskip

\noindent{\bf Acknowledgements.}  The author warmly thank F. Bolley and A. Guillin for stimulating  discussions. This research was supported  by the French ANR-12-BS01-0019 STAB project. The author would like also to thanks the anonymous referee for carefully reading the manuscript.

\newcommand{\etalchar}[1]{$^{#1}$}

\end{document}